\documentclass{elsarticle}

\usepackage[utf8]{inputenc}
\usepackage[english]{babel}

\usepackage{amsmath}
\usepackage{bm}
\usepackage{tikz}
\usepackage{pgfplots}

\usepackage{hyperref}

\newtheorem{thm}{Theorem}
\newtheorem{lem}[thm]{Lemma}
\newdefinition{rmk}{Remark}
\newproof{pf}{Proof}
  
\journal{arXiv.org}

\begin{document}

\begin{frontmatter}

\title{Splitting schemes for unsteady problems involving the grad-div operator \footnotemark[1]}

\author[ua]{Peter Minev\corref{cor}}
\ead{minev@ualberta.ca}

\author[nsi,univ]{Petr N. Vabishchevich}
\ead{vabishchevich@gmail.com}

\address[ua]{Mathematical and Statistical Sciences, 677 Central Academic Building, University of Alberta, Edmonton, Alberta, Canada}
\address[nsi]{Nuclear Safety Institute, Russian Academy of Sciences, 52, B. Tulskaya, Moscow, Russia}
\address[univ]{North-Eastern Federal University, 58, Belinskogo, Yakutsk, Russia}

\cortext[cor]{Corresponding author}

\begin{abstract}
In this paper we consider various splitting schemes for unsteady problems containing the grad-div operator.
The fully implicit discretization of such problems would yield at each time step a linear problem that couples all components of the
solution vector.  In this paper we discuss various possibilities to decouple the equations for the different components
that result in unconditionally stable schemes.  If the spatial discretization uses Cartesian grids, the resulting
schemes are Locally One Dimensional (LOD).  The stability analysis of these schemes is based on the general stability
theory of additive operator-difference schemes developed by Samarskii and his collaborators.  The results of the theoretical
analysis are illustrated on a 2D numerical example with a smooth manufactured solution.
\end{abstract}


\end{frontmatter}

\section{Introduction}

The grad-div ($\nabla \nabla\!\cdot$) operator appears in various problems in science and engineering, 
the obvious examples being the
Navier-Stokes equations in the so-called stress-divergence form (see e.g. \cite{gresho1998incompressible}), 
the equations of linear elasticity, some formulations of the 
Maxwell equations, etc. (see e.g. \cite{LandauLifshic1986,LandauLifshic1982}).
For a more detailed list of applications containing this operator the reader is referred to
\cite{AADM07}.  Sometimes, this operator may appear as a result of a numerical regularization 
of problems involving incompressible fields as in the case of the so-called artificial compressibility methods
 (see for example \cite{Shen_1995} , \cite{GM14}, and \cite{GM16}).
 
 The main issue with the presence of a $\nabla \nabla\!\cdot$ operator is that in implicit discretizations it couples 
 the equations for the various components of the vector field, similarly to the rot-rot ($\nabla\!\times \nabla\!\times$) operator.
 If this coupling can be avoided, the resulting linear system would be easier to solve, no matter if direct or iterative linear solvers are used.  
 Such decoupling of the vectorial problem has already been proposed, for the case of the Maxwell equations involving the $\nabla\!\times \nabla\!\times$ operator, in
 \cite{V04} and \cite{Vab14}, section 3.4.  It results in a set of uncoupled elliptic problems for each component of the solution vector.
 In the present paper we show that the same approach
 works in the case of problems involving $\nabla \nabla\!\cdot$ operators.
If the spatial discretization is performed on Cartesian grids, then these splittings result in Locally One Dimensional (LOD) schemes that are
 very efficient from computational standpoint.  
 Note that while the $\nabla\!\times \nabla\!\times$ operator is positive the $\nabla \nabla\!\cdot$ operator  is only non-negative
 and this difference requires some modifications in the stability analysis of the resulting schemes.
Nevertheless, the stability analysis of the schemes proposed in this paper is based on the general stability theory of additive
operator-difference schemes developed by Samarskii and his collaborators (see \cite{Samarskii1989,SamarskiiMatusVabischevich2002}).
Taking into account the non-negativity of the $\nabla \nabla\!\cdot$ operator, we prove the unconditional stability and provide
{\em a priori} estimates for several decoupling schemes that have been used for vectorial problems with positive operators (see \cite{V04} and \cite{Vab14}, section 3.4).

 The rest of the paper is organized as follows. In the next section we formulate the problem and derive some {\em a priori} estimates for it.
 In section \ref{sec2} we consider some standard two-level schemes for unsteady problems with grad-div operators, and discuss their stability. 
 In section \ref{sec3} two-level splitting schemes with implicit block-diagonal or block-triangular structure are considered, as
 well as alternating block-triangular two- and three-level schemes.  The theoretical results are verified on a numerical example with a manufactured 
 solution in section \ref{sec4}.  Finally, in the last section we summarize the results of this paper.

\section{Problem formulation} \label{sec1}
Consider the problem: Find $\bm u(\bm x,t) = (u_1, ..., u_d)^T$ that satisfies:

\begin{subequations} 
\begin{align}
   \frac{\partial \bm u}{\partial t} - {\rm grad} ( k({\bm x}) \, {\rm div}  \, \bm u ) &= f({\bm x}, t) , 
   \quad \bm x \in \Omega , \quad 0 < t \leq T,\label{1}\\
   (\bm u \cdot \bm n) &= 0, \quad {\bm x} \in \partial\Omega, \label{2}\\
   \bm u({\bm x}, 0) &= \bm u_0({\bm x}),
  \quad {\bm x} \in \Omega, \label{3}
  \end{align}
  \label{123}
\end{subequations}

\noindent where $\Omega \subset R^d$, $d=2,3$ is a bounded polygonal domain with a Lipschitz continuous boundary 
$\partial\Omega$, and $\bm n$ is the outward normal to $\partial \Omega$. The coefficient $k({\bm x}) \geq  0$ 
and the source term $f({\bm x}, t)$ are supposed to be sufficiently smooth.

Let $(\cdot,\cdot)$ denote the standard $L^2$ inner product over $\Omega$, 
and $\|\cdot\|$ be the corresponding norm
for scalar and vector functions $v(\bm x)$ and $\bm v(\bm x)$, correspondingly:
\[
 (v, w) = \int_{\Omega } v(\bm x)  w(\bm x) d \bm x,
 \quad \|v\| = (v,v)^{1/2} ,
\] 
\[
 (\bm v, \bm w) = \sum_{i=1}^{d} (v_i, w_i) ,
 \quad \|\bm v\| = \Big (\sum_{i=1}^{d}\|v_i\|^2 \Big )^{1/2} .
\] 
It is quite obvious that problem (\ref{123})  is well posed and we can easily obtain {\em a priori} estimates for it.
To reduce the complexity of the notation in the  paper we use calligraphic letters for denoting operators in infinite 
dimensional spaces and standard capital letters for their finite dimensional approximations.
Equation (\ref{1}) on the set of functions satisfying (\ref{2}) can be reformulated as a standard Cauchy problem
for a first order evolutionary equation:
\begin{subequations}\label{45}
\begin{align}
 \frac{d \bm u}{d t} + \mathcal{A} \bm u &= \bm f(t),  \quad 0 < t \leq T, \label{4}\\
 \bm u(0) &= \bm u_0 , \label{5}
  \end{align}
\end{subequations} 
where $\bm u(t) = \bm u(\cdot,t)$.
It is clear that  $\mathcal{A}$ is self-adjoint and non-negative in $L^2(\Omega)$ i.e.:
\begin{equation}\label{6}
 \mathcal{A} = \mathcal{A}^* \geq 0 .
\end{equation} 

To obtain an {\em a priori} estimate for the solution of (\ref{45}), (\ref{6}),  we multiply (\ref{4}) by $\bm u$ to obtain:
\[
 \left (\frac{d \bm u}{d t}, \bm u \right ) + (\mathcal{A} \bm u, \bm u)  = (\bm f,  \bm u) .
\]  
Taking into account (\ref{6}) and the identity 
\[
 \left (\frac{d \bm u}{d t}, \bm u \right ) = \frac{1}{2}  \frac{d }{d t} \|\bm u\|^2 =
  \|\bm u\| \frac{d }{d t} \|\bm u\| 
\] 
we obtain that:
\[
 \frac{d }{d t} \|\bm u\| \leq  \|\bm f\| .
\] 
An application of the Gronwall inequality yields the stability estimate for problem (\ref{45}), (\ref{6}):
\begin{equation}\label{7}
 \|\bm u(t)\| \leq  \| \bm u_0\| + \int\limits_{0}^{t}  \|\bm f(\theta)\| d\theta. 
\end{equation} 
In addition to the estimate (\ref{7}), we can also derive an estimate that is suitable for more general equations
by multipying (\ref{4}) by $\displaystyle \frac{d \bm u}{d t}$ to obtain:
\begin{equation}\label{8}
 \left (\frac{d \bm u}{d t}, \frac{d \bm u}{d t} \right ) + 
 \frac{1}{2} \frac{d }{d t} | \bm u|_{\mathcal{A}}^2 = \left (\bm f,  \frac{d \bm u}{d t} \right ) , 
\end{equation} 
where $| \cdot |_{\mathcal{A}}$ is the seminorm induced by $\mathcal{A}$:
$| \bm u|_{\mathcal{A}}^2 = (\mathcal{A} \bm u, \bm u)$.
Estimating the right hand side of (\ref{8}) by the inequality:
\[
 \left (\bm f,  \frac{d \bm u}{d t} \right ) \leq \frac{1}{2} \|\bm f\|^2
 + \frac{1}{2} \left (\frac{d \bm u}{d t}, \frac{d \bm u}{d t} \right ) , 
\] 
substituting it in (\ref{8}), and integrating for  $0<t \leq T$, we readily obtain the estimate:
\begin{equation}\label{9}
 \int_{0}^{T} \left \| \frac{d \bm u}{d t} \right \|^2 d t 
 + |\bm u(T)|_{\mathcal{A}}^2 \leq |\bm u_0|_{\mathcal{A}}^2
 + \int_{0}^{T}  \|\bm f(t) \|^2 d t . 
\end{equation} 
It is the discrete version of (\ref{9}), that will guarantee the stability of the schemes considered in the next two sections.

\begin{rmk}\label{rm1}
From (\ref{9}) we can obtain a uniform in $t$ estimate of type (\ref{7}), for the square of the norm of the solution.
Indeed, for a differentiable function $g(t)$ we have that:
\[
 g(t) = \int_{0}^{t} \frac{d g(\theta)}{d \theta } d \theta + g(0) ,
\] 
which implies that:
\[
 g^2(t) \leq  2 \left ( \int_{0}^{t} \frac{d g(\theta)}{d \theta } d \theta \right )^2 + 2 g^2(0) .
\] 
Taking into account the inequality:
\[
 \left ( \int_{0}^{t} \frac{d g(\theta)}{d \theta } d \theta \right )^2 \leq 
 t \int_{0}^{t} \left (  \frac{d g(\theta)}{d \theta } \right )^2 d \theta  ,
\] 
we readily obtain that:
\begin{equation}\label{10}
 g^2(t) \leq  2 t \int_{0}^{t} \left (  \frac{d g(\theta)}{d \theta } \right )^2 d \theta + 2 g^2(0) . 
\end{equation} 
Using the estimates (\ref{10}) with $g=\bm u$, and (\ref{9}), we obtain the following estimate for the square of the $L^2$-norm of the solution:
\begin{equation}\label{11}
  \|\bm u(t)\|^2 \leq 2 \|\bm u_0\|^2 + 2 t |\bm u_0|_{\mathcal{A}}^2 + 2 t \int\limits_{0}^{t}  \|\bm f(\theta)\|^2 d\theta .
\end{equation} 
\end{rmk}

\section{Standard time discretization} \label{sec2}

The component form of  (\ref{4}) is given by:
\begin{equation}\label{12}
 \frac{d u_i}{d t} + \sum_{j=1}^{d} \mathcal{A}_{ij} u_j = f_i(t),
 \quad i = 1,..., d , 
 \quad 0 < t \leq T, 
\end{equation}   
where:
\[
 \mathcal{A}_{ij} u_j = - \frac{\partial }{\partial x_i} \left (
 k(\bm x) \frac{\partial u_j}{\partial x_j} \right ),
 \quad i,j = 1,...., d .   
\] 
For the operator $\mathcal{A}$ we have that
\[
 \sum_{i,j=1}^{d} (\mathcal{A}_{ij} u_j, u_i) =
 \sum_{i,j=1}^{d} \left ( k(\bm x) \frac{\partial u_j}{\partial x_j} ,  \frac{\partial u_i}{\partial x_i} \right ) =
 \sum_{i=1}^{d}  \left ( k(\bm x) \left ( \frac{\partial u_i}{\partial x_i} \right )^2 , 1 \right ) \geq 0, 
\] 
i.e. it satisfies conditions (\ref{6}).  $\mathcal{A}$ can also be bounded from above by its block-diagonal part as follows:
\begin{equation}\label{13}
 2 (\mathcal{A}_{ij} u_j, u_i) \leq (\mathcal{A}_{ii} u_i, u_i) +
  (\mathcal{A}_{jj} u_j, u_j) ,
 \quad i,j = 1, ..., d , 
\end{equation} 
i.e.:
\[
 \sum_{i,j=1}^{d} (\mathcal{A}_{ij} u_j, u_i)
 \leq d \sum_{i=1}^{d} (\mathcal{A}_{ii} u_i, u_i) .
\]
This estimate guarantees that:
\begin{equation}\label{14}
 \mathcal{A} \leq d \, \mathcal{D} ,
\end{equation} 
where $\mathcal{D} = \mathrm{diag} ( \mathcal{A}_{11}, ...,  \mathcal{A}_{dd})$.

Problem (\ref{45}) can be discretized in space by means of standard $H^{div}$-stable finite elements (Raviart-Thomas, BDM), finite volume 
or standard finite difference methods, and we write the Cauchy problem for the resulting ODE system as:
\begin{subequations}\label{1516}
\begin{align}
 \frac{d \bm v}{d t} + A \bm v &= \bm \varphi (t), \quad 0 < t \leq T, \label{15}\\
 \bm v(0) &= \bm v_0 . \label{16}
 \end{align}
\end{subequations} 
Note that in case of non-uniform grids, the time derivative in (\ref{15}) would be multiplied by a mass matrix, approximating the identity operator.  
The appearance of such a matrix almost never alters the considerations provided further in the paper, and therefore, we mostly omit it. However, whenever
needed, we comment on the role of the mass matrix.
Besides, we assume that the operator ${A}=(A_{ij})$  inherits the main properties of $\mathcal{A}$, (\ref{6})  and (\ref{14}),  i.e.:
\begin{equation}\label{17}
 A = A^* \geq 0,
\end{equation} 
\begin{equation}\label{18}
 A \leq d \, D,
 \quad D = \mathrm{diag} ( A_{11}, ...,  A_{dd}) ,
\end{equation} 
in the corresponding finite dimensional space.

We first consider a two-level approximation of problem (\ref{1516}) (see, e.g. \cite{Samarskii1989}).
Let $\tau$ be a step-size of a uniform grid in time such that 
$\bm v^n = \bm v(t^n), \ t^n = n \tau$, $n = 0,1, ..., N, \ N\tau = T$.
The general form of the two-level scheme is given by:
\begin{subequations}\label{1920}
\begin{align}
 \frac{\bm v^{n+1} - \bm v^{n}}{\tau} + 
 A (\sigma \bm v^{n+1} + (1-\sigma) \bm v^n ) &= 
 \bm \varphi (\sigma t^{n+1} + (1-\sigma) t^n) ,\label{19}\\ 
  \bm v^0 &= \bm v_0 ,\label{20}
 \end{align}
\end{subequations} 
where $\sigma $ is the weight of the scheme.
Obviously, if $\sigma = 0, 0.5, 1$, the scheme turns into the  Euler explicit, Crank--Nicolson,
and  Euler implicit schemes, respectively. 

In order to study the stability of (\ref{1920}) and the splitting schemes that follow below, we will use the general approach for two- and three-level schemes developed in \cite{Samarskii1968,Samarskii1970}.
If the operator $A$ is positive definite then we can use the sufficient and necessary conditions for stability obtained in \cite{Samarskii1989,SamarskiiMatusVabischevich2002}.  However,
in the present case it is only semidefinite and we cannot apply this condition.  Nevertheless, as demonstrated below, a modification of this analysis easily yields a sufficient condition for stability.
To begin, we notice that  (\ref{1516}), (\ref{17}), as well as all two-level splitting schemes discussed further in the paper, can be written in the following canonical form:
\begin{equation}\label{21}
 B \frac{\bm v^{n+1} - \bm v^{n}}{\tau} + 
 A  \bm v^n =  \bm \psi^n , 
 \quad n = 0,1, ..., N-1 , 
\end{equation} 
subject to the initial condition (\ref{20}). Here the operator $A$ is assumed to be positive semi-definite and $B$ is positive definite but in some of the splitting schemes below it is not self-adjoint.  
This is why we need to consider the self-adjoint part of it given by:
\[
 B_0 = \frac{1}{2} (B + B^*) .
\] 
Further, with each positive definite and self-adjoint operator $C$  we associate a Hilbert space $H_C$, with an inner product and norm defined by:
\[
 (\bm v,\bm w)_C = (C \bm v,\bm w), 
 \quad \|\bm v\|_C = (C\bm v, \bm v)^{1/2} . 
\] 
Then, a sufficient condition for the stability of the general two-level scheme given by (\ref{21}), (\ref{20}), is given by the following lemma:

\begin{lem}\label{l-1}
Let in (\ref{21})
\begin{equation}\label{22}
 C = B_0 - \frac{\tau }{2} A > 0,
\end{equation} 
then the scheme (\ref{21}), (\ref{20}) is stable in $H_C$, and we have the following {\em a priori} estimate:
\begin{equation}\label{23}
 \sum_{n=0}^{N-1} \tau  \left \| \frac{\bm v^{n+1} - \bm v^{n}}{\tau} \right \|^2_C + |\bm v^{N}|^2_A 
 \leq |\bm v^{0}|^2_A + \sum_{n=0}^{N-1} \tau \|\bm \psi^n\|^2_{C^{-1}} .   
\end{equation}  
\end{lem}

\begin{pf}
We first rewrite (\ref{21}) as:
\[
 \left ( B - \frac{\tau }{2} A \right ) \frac{\bm v^{n+1} - \bm v^{n}}{\tau} + 
 A  \frac{\bm v^{n+1} + \bm v^{n}}{2} =  \bm \psi^n .
\] 
Multiplying it by $2(\bm v^{n+1} - \bm v^{n})$, and using (\ref{22}) we obtain:
\begin{equation}\label{24}
 2 \tau  \left \| \frac{\bm v^{n+1} - \bm v^{n}}{\tau} \right \|^2_C + 
 |\bm v^{n+1}|^2_A -  |\bm v^{n}|^2_A = 
 2 \tau \left (\bm \psi^n , \frac{\bm v^{n+1} - \bm v^{n}}{\tau} \right ) .  
\end{equation} 
The right hand side of this equation is estimated as follows:
\[
 2 \left (\bm \psi^n , \frac{\bm v^{n+1} - \bm v^{n}}{\tau} \right ) 
 \leq  \left \| \frac{\bm v^{n+1} - \bm v^{n}}{\tau} \right \|^2_C
 + \|\bm \psi^n\|^2_{C^{-1}} ,
\] 
thus yielding the estimate:
\begin{equation}\label{24.1}
 \tau  \left \| \frac{\bm v^{n+1} - \bm v^{n}}{\tau} \right \|^2_C + 
 |\bm v^{n+1}|^2_A -  |\bm v^{n}|^2_A \leq  \tau \|\bm \psi^n\|^2_{C^{-1}}  .  
\end{equation} 
Summing  (\ref{24.1}) for $n=0,1, ..., N-1$ yields the requaired estimate (\ref{23}).
\end{pf}

Using lemma \ref{l-1}, we can easily study the stability of the scheme (\ref{1920}).
In this case we clearly have that:
\[
 B = I + \sigma \tau A ,
\] 
and
\begin{equation}\label{25}
 C = I + \left ( \sigma - \frac{1}{2} \right ) \tau A ,
\end{equation} 
 with $I$ being the identity operator.  Then an application of lemma \ref{l-1}
 easily proves the following theorem.
 
\begin{thm}\label{t-1}
Let in (\ref{19}) $\sigma \geq 1/2$,
then the scheme (\ref{1920}) is stable in $H_C$, where $C$ is given by (\ref{25}).
The solution satisfies the a priori estimate (\ref{23}), where:
\[
 \bm \psi^n = \bm \varphi (\sigma t^{n+1} + (1-\sigma) t^n) , 
 \quad n = 0,1, ..., N-1 .  
\] 
\end{thm}

Note that the estimate (\ref{23}) is the discrete version of  (\ref{9}).

\begin{rmk}
Repeating the argument from remark \ref{rm1} in the discrete setting, we can obtain a discrete version of the estimate (\ref{11}) for the scheme  (\ref{21}), (\ref{20}).
Indeed, if we introduce the grid function $g^n, \ n = 0,1, ..., N$, we have that:
\[
 g^{n+1} = \sum_{k=0}^{n} \tau \frac{g^{k+1} - g^k}{\tau} + g^0
\]  
and therefore:
\[
 (g^{n+1})^2 \leq 2\left ( \sum_{k=0}^{n} \tau \frac{g^{k+1} - g^k}{\tau} \right )^2  + 2 (g^0)^2 .
\]  
Using the Cauchy--Schwarz inequality
\[
 \left ( \sum_{k=0}^{n} a_k b_k \right )^2 
 \leq  \sum_{k=0}^{n} a_k^2 \sum_{k=0}^{n} b_k^2 
\] 
with
\[
  a_k = \tau^{1/2}, 
  \quad   b_k = \tau^{1/2} \frac{g^{k+1} - g^k}{\tau} ,
\] 
we obtain the discrete analog of (\ref{10}) 
\begin{equation}\label{26}
 (g^{n+1})^2 \leq 2t^{n+1} \sum_{k=0}^{n} \tau \left ( \frac{g^{k+1} - g^k}{\tau} \right )^2  + 2 (g^0)^2. 
\end{equation} 
Substituting it, with $g={\bm v}$, into (\ref{24.1}) gives:
\begin{equation}\label{27}
 \|\bm v^{n+1}\|^2_C  \leq 2 \|\bm v^{0}\|^2_C + 2t^{n+1} |\bm v^{0}|^2_A +
 2t^{n+1} \sum_{k=0}^{n} \tau \|\bm \psi^k\|^2_{C^{-1}} ,
\end{equation} 
which is the discrete version of the estimate (\ref{11}).
\end{rmk}

\section{Splitting schemes} \label{sec3}
In order to obtain the solution to (\ref{1920}), at each time level $n+1$, we need to solve linear systems of type:
\begin{equation}\label{28}
 v_i^{n+1} + \tau \sigma \sum_{j=1}^{d} A_{ij} v_j^{n+1} = r_i^{n},
 \quad i = 1,..., d ,  
\end{equation} 
that, if $\sigma >0$, couple implicitly all components  $v_i^{n+1}, \ i = 1,..., d$.
Splitting schemes are aimed at decoupling of the problems for the different components by approximating some of the 
terms in the right hand side of (\ref{28}) explicitly (see \cite{Vab14}) .
Then, instead of (\ref{28}) we solve decoupled problems of type:
\begin{equation}\label{29}
 v_i^{n+1} + \tau \sigma_1 \sum_{j=1}^{i-1} A_{ij} v_j^{n+1} + \tau \sigma_2 A_{ii} v_i^{n+1}= r_i^{n},
 \quad i = 1,..., d ,  
\end{equation} 
with $\sigma_1$ and  $\sigma_2$ being some consistently chosen weights.
In particular, if $\sigma_1 = 0$ the problems for all $v_i^{n+1}, \ i = 1,..., d$,
can be solved independently and this facilitates the parallel implementation of the algorithm.  It is also clear that
the problems (\ref{29}) require the solution of one dimensional problems only and therefore the resulting schemes are LOD.

\subsection{Block Jacobi splitting scheme} 

In this case all terms in (\ref{1516}), coming from the off-diagonal blocks of the operator
$A$ are discretized explicitly to obtain the scheme:
\begin{equation}\label{30}
 \frac{\bm v^{n+1} - \bm v^{n}}{\tau} + 
 D (\sigma \bm v^{n+1} + (1-\sigma) \bm v^n) + (A-D) \bm v^n = 
 \bm \psi^n .
\end{equation} 
It is clear that in this case we need to solve problems of type (\ref{29}) with $\sigma_1 = 0$ and $\sigma_2 = \sigma$.

\begin{thm}\label{t-2}
If $\sigma \geq d/2$ then the scheme (\ref{30}), (\ref{20}), is stable in $H_C$  with:
\begin{equation}\label{31}
  C = I +  \sigma \tau D - \frac{\tau}{2}  A .
\end{equation} 
Moreover, its solution obeys the estimate (\ref{23}).  
\end{thm}

\begin{pf}
Obviously equation (\ref{30}) can be rewritten in the form of (\ref{21}) with:
\[
 B = I + \sigma \tau D .
\] 
Then taking into account the inequalities (\ref{22}) with $B_0=B$, and (\ref{18}), and the definition of $C$ given by (\ref{31}) 
we readily obtain that:
\[
 C = C^* \geq I
\] 
for $\sigma \geq d/2$. Then lemma  \ref{l-1} immediately yields the estimate (\ref{23}) and therefore the stability of the block-Jacobi scheme.  
\end{pf}

\subsection{Block Gauss--Seidel splitting scheme} 

In this section we consider the possibility to create a stable scheme that leads to linear problems of type
(\ref{29}) for $\sigma_1 > 0$. For this, we first decompose $A$ into a block lower-triangular, diagonal, and upper-triangular parts as follows:
\begin{equation}\label{31.1}
 A = L + D + U,
 \quad L^* = U, 
\end{equation}
where
\[
 L = (L_{ij}),
 \quad  L_{ij} = \left \{
\begin{array}{cc}
  A_{ij},   &  j < i \\
  0,  & j \geq i ,
\end{array}
\quad i, j = 1,..., d ,
\right .  
\] 
and
\[
 U = (U_{ij}),
 \quad  U_{ij} = \left \{
\begin{array}{cc}
   0,   &  j \leq  i \\
 A_{ij},  & j > i ,
\end{array}
\quad i, j = 1,..., d .
\right .  
\] 
Then the block Gauss--Seidel type scheme for (\ref{1516}) is given by:
\begin{equation}\label{32}
 \frac{\bm v^{n+1} - \bm v^{n}}{\tau} + 
 (L + D) \bm v^{n+1} + U \bm v^n = 
 \bm \psi^n .
\end{equation} 
Its stability is guaranteed by the following theorem.

\begin{thm}\label{t-3}
The scheme (\ref{32}), (\ref{20}) is unconditionally stable in $H_C$ with
\begin{equation}\label{33}
  C = I + \frac{\tau}{2}  D ,
\end{equation} 
and its solution satisfies the estimate (\ref{23}).  
\end{thm}

\begin{pf}
The result again follows immediately from lemma \ref{l-1}.
We first rewrite (\ref{32}) into the two-level canonical form (\ref{21}) with
\[
 B = I + \tau (L + D) .
\] 
Then clearly
\[
 B_0 = \frac{1}{2}\left(B+B^*\right) = I + \frac{\tau}{2}  D + \frac{\tau}{2}  A, 
\] 
and therefore $\displaystyle C=B_0-\frac{\tau}{2} A=I + \frac{\tau}{2}  D$.  This is obviously a positive definite operator 
and we can imply the conclusion of lemma  \ref{l-1}, hence concluding the proof.
\end{pf}

\subsection{Alternating triangular splitting scheme} 

This scheme consists of two block--Gauss--Seidel steps, the first one using a block lower-triangular matrix $B$ and the second one
using its transpose. This choice justifies the name of the scheme.  It in fact uses the same idea as the one that leads to explicit 
unconditionally stable schemes for second order parabolic equations proposed by \cite{Saulev1960,Samarskii1964}. However, in the
present case we have to use the block triangular decomposition of $A$ and therefore the scheme is not explicit.
Using the decomposition given by (\ref{31.1}) we first split $A$ as:
\begin{equation}\label{34}
 A = A_1 + A_2,
 \quad A_1^* = A_2, 
\end{equation} 
where
\[
 A_1 = L + \frac{1}{2} D,
 \quad  A_2 = U + \frac{1}{2} D .
\] 
Then the alternating triangular scheme is produced from (\ref{21})  using an operator $B$ given by:
\begin{equation}\label{35}
 B = (I + \tau \sigma A_1) (I + \tau \sigma A_2) .
\end{equation} 
Its stability is analized by the following theorem.

\begin{thm}\label{t-4}
If $\sigma \geq 1/2$, the scheme comprised by  (\ref{21}), (\ref{20}), (\ref{34}), and (\ref{35}) is stable in $H_C$ with
\begin{equation}\label{36}
 C = I + \left ( \sigma - \frac{1}{2} \right ) \tau A + \sigma^2 \tau^2 A_1 A_2 ,
\end{equation} 
and its solution satisfies the estimate (\ref{23}).  
\end{thm}
\begin{pf}
The proof follows along the same lines as in theorem \ref{t-1}, since,
using the second identity in (\ref{34}), it is straightforward to show that $C = C^* > 0$.
\end{pf}

An advantage of the alternating triangular splitting scheme 
is that the choice $\sigma = 1/2$ yields a second order unconditionally stable scheme. 

\begin{rmk}
Note that if the spatial discretization uses non-uniform Cartesian grids or finite elements,
the identity operator would be approximated by a mass matrix $M$ which is not equal to $I$.
This change requires that we modify the choice of $B$ in such case to:
\[
 B = (M + \tau \sigma A_1)M^{-1} (M + \tau \sigma A_2) .
\] 
\end{rmk}
However,  the stability result in  theorem \ref{t-4} remains valid if we modify $C$ in (\ref{36}) to
$ C = M + \left ( \sigma - \frac{1}{2} \right ) \tau A + \sigma^2 \tau^2 A_1 M^{-1} A_2$.

\subsection{Three-level alternating triangular scheme} 
 
Recently, a three-level modification of the alternating triangular method has been proposed in \cite{vab2014}.
Its consistency differs from the consistency  of the standard scheme  (\ref{19}) by a term that is of order of $\mathcal{O}(\tau^3)$
i.e. its splitting error is third order. 
Such a three-level alternating triangular scheme in the present case can be written as:
\begin{equation}\label{37}
\begin{split}
 (I + \tau \sigma A_1) (I + \tau \sigma A_2) \frac{\bm v^{n+1} - \bm v^{n}}{\tau} -
 & \sigma^2 \tau^2 A_1 A_2  \frac{\bm v^{n} - \bm v^{n-1}}{\tau} + 
 A  \bm v^n =  \bm \psi^n , \\
 &  n = 0,1, ..., N-1 , 
\end{split} 
\end{equation} 
where $\bm v^{-1}, \bm v^0$ must be given (usually produced by a two-level scheme).
Note that it can also be written in the following canonical three-level form (see \cite{Vab14}, section 2.3):
\begin{equation}\label{38}
 B \frac{\bm v^{n+1} - \bm v^{n-1}}{2\tau} +
 R \frac{\bm v^{n+1} - 2\bm v^{n} + \bm v^{n-1}}{\tau} + 
 A  \bm v^n =  \bm \psi^n ,
\end{equation} 
where 
\[
 B = I + \tau \sigma A,
 \quad R = \frac{1}{2} (I + \tau \sigma A) + \sigma^2 \tau^2 A_1 A_2 .
\] 
Also $R^* = R > 0$. 
Using this form, it is possible to apply some of the general results in \cite{Samarskii1989,SamarskiiMatusVabischevich2002}  
of the analysis of three-level schemes in such a canonical form,  
and obtain {\em a priori} estimates in the present case.  However, the norms used in these estimates would be quite complicated.  
Here we proceed somewhat differently in this particular case, similarly to the approach in lemma \ref{l-1}.  
We first re-write (\ref{37})  as follows:
\begin{equation}\label{39}
 C \frac{\bm v^{n+1} - \bm v^{n}}{\tau} + 
 A  \frac{\bm v^{n+1} + \bm v^{n}}{2} +
 R \frac{\bm v^{n+1} - 2\bm v^{n} + \bm v^{n-1}}{\tau} = \bm \psi^n ,
\end{equation} 
where $\sigma \geq 1/2$, and $C$ is the same as in (\ref{25}).
In order to establish an {\em a priori} estimate for (\ref{39}) we multiply it by  $2 \tau \bm w^{n+1}$, where $\displaystyle 
 \bm w^{n+1} = \frac{\bm v^{n+1} - \bm v^{n}}{\tau}$, and similarly to (\ref{24}) we obtain:
\begin{equation}\label{40}
 2 \tau  \| \bm w^{n+1} \|^2_C + 
 |\bm v^{n+1}|^2_A -  |\bm v^{n}|^2_A +
 2 \tau (R (\bm w^{n+1}-\bm w^{n}), \bm w^{n+1}) =  
 2 \tau (\bm \psi^n , \bm w^{n+1} ) .  
\end{equation} 
The right hand side can be estimated as:
\[
 2 (\bm \psi^n , \bm w^{n+1} ) \leq \|\bm w^{n+1} \|_C^2 + \| \bm \psi^n \|^2_{C^{-1}} .
\] 
Taking into account that
\[
 2 (R (\bm w^{n+1}-\bm w^{n}), \bm w^{n+1}) = 2 \|\bm w^{n+1}\|_R^2 - 2(R\bm w^{n}, \bm w^{n+1}) \geq 
 \|\bm w^{n+1}\|_R^2 - \|\bm w^{n}\|_R^2 ,
\] 
and substituting this result in  (\ref{40}) we readily obtain:
\[
 \tau  \| \bm w^{n+1} \|^2_C + |\bm v^{n+1}|^2_A + \tau  \|\bm w^{n+1}\|_R^2 \leq  
 |\bm v^{n}|^2_A +  \tau  \|\bm w^{n}\|_R^2 + \tau \| \bm \psi^n \|^2_{C^{-1}} .
\] 
Summing these inequalities for $n = 0,1, ..., N-1$ gives:
\[
 \sum_{n=0}^{N-1} \tau  \| \bm w^{n+1} \|^2_C + |\bm v^{N}|^2_A +  \tau  \|\bm w^{N}\|_R^2 \leq 
 |\bm v^{0}|^2_A +  \tau  \|\bm w^{0}\|_R^2 + \sum_{n=0}^{N-1} \tau \| \bm \psi^n \|^2_{C^{-1}} .
\] 
Substituting $\displaystyle \bm w^{n+1} = \frac{\bm v^{n+1} - \bm v^{n}}{\tau}$, we finally obtain the estimate:
\begin{equation}\label{41}
\begin{split}
 \sum_{n=0}^{N-1} \tau  \left \| \frac{\bm v^{n+1} - \bm v^{n}}{\tau} \right \|^2_C  & + |\bm v^{N}|^2_A +
 \sigma^2 \tau^3 \left \| A_2 \frac{\bm v^{N} - \bm v^{N-1}}{\tau} \right \|^2 \\
 & \leq |\bm v^{0}|^2_A + 
 \sigma^2 \tau^3 \left \| A_2 \frac{\bm v^{0} - \bm v^{-1}}{\tau} \right \|^2 +
 \sum_{n=0}^{N-1} \tau \|\bm \psi^n\|^2_{C^{-1}} .  
\end{split}
\end{equation} 
This result is the three-level counterpart of the estimate (\ref{23}) valid for the two-level schemes discussed above, and
can be summarized in the following theorem:
\begin{thm}\label{t-5}
If $\sigma \geq 1/2$, then the scheme (\ref{37}) is stable in  $H_C$, where $C$ is given by  (\ref{25}).
More precisely its solution satisfies the estimate in (\ref{41}).
\end{thm}

\section{Numerical results}\label{sec4}

We test all schemes considered in the previous section on a manufactured solution 
of (\ref{123})   given by: 
\[
 \bm u = (\sin t \sin(\pi x_1) \sin(\pi x_2), \cos t \sin(\pi x_1) \sin(\pi x_2))^T 
\] 
in the domain $\Omega = (0,1)\times (0,1)$, and with a properly chosen source term. 
We use a discretization based on the so called MAC stencil i.e. the grid $\omega$ for the the
two components of $\bm u$ is staggered in each direction. 
The grid size is constant and the solution that is used to measure the error is interpolated in the centroids of the MAC 
cells.  The error is measured in the $L_2(\omega)$ norm which, for a grid function $\phi(\bm x_\alpha), \bm x_\alpha \in \omega$, 
is given by: 
\[
 \|\phi\|^2 = \sum \limits_{\bm x_\alpha \in \omega} \phi^2(\bm x_\alpha) h^2 ,
\]
where $h$ is the grid size.
For reference we also provide results with the Crank-Nicolson discretization given by (\ref{19}) with $\sigma =1/2$.  

The results for the error of the various approximations at $t=10$ computed
on a grid of $200 \times 200$ MAC cells are presented in table \ref{tb:1}.
They clearly demonstrate that: (i) the only scheme that yields relatively good results, as compared to the Crank-Nicolson scheme,
is the three-level scheme; (ii) as expected, the Jacobi and Gauss-Seidel schemes have almost first order of convergence; (iii) the 
alternating triangular scheme
has a better-than-first order of convergence; (iv) the three-level scheme has a second order of convergence.  It must be pointed out
that the three-level scheme requires the solution of 1D parabolic problems only, no matter what is the dimensionality of the problem,
while the Crank-Nicolson scheme would be much more computationally intensive.

\begin{table}
  \caption{Discrete norm of the error for $\bm u$ at $T=10$, $\nu=1$, and a grid size of $0.005$. 
J ($\sigma=1$) stands for the Jacobi scheme 
  with $\sigma = 1$, GS ($\sigma=1$) stands for the Gauss-Seidel scheme, AT ($\sigma=0.5$) stands for the alternating triangular scheme, 
  3-level ($\sigma=0.5$) stands for the three-level scheme, and CN stands for the Crank-Nicolson scheme.}
\centering
\begin{tabular}{|c|c|c|c|c|c|} 
\hline
$\tau $ & J ($\sigma=1$)  & GS ($\sigma=1$) & AT ($\sigma=0.5$) & 3-level ($\sigma=0.5$) & CN  \\
\hline
0.1 & 0.144 & 0.2 & 0.185 & 0.01 & 6 $\times 10^{-4}$  \\ \hline
0.05 & 0.103 & 0.129 & 0.067 & 1.9$\times 10^{-3}$  & 1.6 $\times 10^{-4}$  \\  \hline
0.025 &  0.066 & 0.075 & 0.023 &   3$\times 10^{-4}$  & 4.6$ \times 10^{-5}$  \\  \hline
0.0125 & 0.039 &  0.041 &  0.008 & 5.7$ \times 10^{-5}$  & 1.9 $\times 10^{-5}$ \\ \hline
0.00625 & 0.021 & 0.021 & 2.7$\times 10^{-3}$ &  1.4 $ \times 10^{-5}$ & 1.2$\times 10^{-5}$ \\ \hline
\end{tabular}
\label{tb:1}
\end{table}

\section{Conclusions}
In this paper we propose several splitting schemes for parabolic problems involving a (time-independent) grad-div operator. 
We proved that all schemes are  stable for a proper choice of the weighting parameter $\sigma$.  
The numerical results, however, clearly demonstrate that the alternating triangular and the three-level  alternating triangular  scheme
converge much faster than the other two schemes.  A clear advantage of these schemes is that their
computational effort only involves the solution of one-dimensional parabolic problems despite the dimension of the original
problem which is significantly less than the effort of an unsplit Crank-Nicolson scheme.

\section*{Acknowledgements}

This work was supported by the Russian Foundation for Basic Research (projects 14-01-00785), by  a Discovery grant of the Natural Sciences and
  Engineering Research Council of Canada, and by a grant \# 55484-ND9 of the Petroleum Research Fund of the American Chemical Society.

\end{document}